\documentclass{amsart}
\usepackage{array}
\usepackage{multirow}
\usepackage{graphicx}
\usepackage{epsfig}
\usepackage{amsfonts,amsmath}
\usepackage{amssymb}
\usepackage{latexsym}
\usepackage{epsf}
\usepackage{graphicx,color,graphics}
\usepackage{amsmath,amssymb}
\usepackage{dsfont}
\usepackage[latin1]{inputenc}  
\newtheorem{theorem}{Theorem}[section]

\numberwithin{equation}{section}

\begin{document}

\title[stability of a weak dissipative Bresse system]{Polynomial and non exponential stability of a weak dissipative Bresse system}
\author[A. Guesmia]{Aissa Guesmia}
\maketitle

\pagenumbering{arabic}

\begin{center}
Institut Elie Cartan de Lorraine, UMR 7502, Universit\'e de Lorraine\\
3 Rue Augustin Fresnel, BP 45112, 57073 Metz Cedex 03, France
\end{center}
\begin{abstract}
In this paper, we study the Bresse system in a bounded domain with linear frictional dissipation working only on the veridical displacement. The longitudinal and shear angle displacements are free. Our first main result is to prove that, independently from the velocities of waves propagations, this linear frictional dissipation does not stabilize exponentially the whole Bresse system. Our second main result is to show that the solution converges to zero at least polynomially. The proof of the well-posedness of our system is based on the semigroup theory. The stability results will be proved using a combination of the energy method and the frequency domain approach.
\end{abstract}

\vskip0,1truecm
Frequency domain approach.
\vskip0,1truecm
{\bf AMS Classification.} 35B40, 35L45, 74H40, 93D20, 93D15.

\renewcommand{\thefootnote}{}
\footnotetext{E-mail addresse: aissa.guesmia@univ-lorraine.fr.}

\section{Introduction}

The subject of this paper is studying the stability of Bresse system under linear frictional damping effective only on the the veridical displacement. This system is defined in $\left( 0,1\right) \times \left(0,\infty \right)$ and takes the form 
\begin{equation}
\left\{
\begin{array}{ll}
\rho _{1}\varphi _{tt}-k\left( \varphi _{x}+\psi +l\,w\right)
_{x}-lk_{0}\left( w_{x}-l\varphi \right) +\delta \varphi _{t}=0 , \vspace{0.2cm}\\
\rho _{2}\psi _{tt}-b\psi _{xx}+k\left( \varphi _{x}+\psi +l\,w\right) =0 , \vspace{0.2cm}\\
\rho _{1}w_{tt}-k_{0}\left( w_{x}-l\varphi \right) _{x}+lk\left( \varphi
_{x}+\psi +l\,w\right) =0 
\end{array}
\right.  \label{syst1}
\end{equation}
along with the initial data
\begin{equation}
\left\{
\begin{array}{ll}
\varphi \left( x,0\right) =\varphi _{0}\left( x\right) ,\,\varphi _{t}\left(
x,0\right) =\varphi _{1}\left( x\right) & \text{in }\left( 0,1\right) , \vspace{0.2cm}\\
\psi \left( x,0\right) =\psi _{0}\left( x\right) ,\,\psi _{t}\left(
x,0\right) =\psi _{1}\left( x\right) & \text{in }\left( 0,1\right) , \vspace{0.2cm}\\
w\left( x,0\right) =w_{0}\left( x\right) ,\,w_{t}\left( x,0\right)
=w_{1}\left( x\right) & \text{in }\left( 0,1\right) 
\end{array}
\right.  \label{cdt_10}
\end{equation}
and the homogeneous Dirichlet-Neumann-Neumann boundary conditions
\begin{equation}
\left\{
\begin{array}{ll}
 \varphi \left( 0,t\right) =\,\psi _{x}\left( 0,t\right) =\,w_{x}\left(
0,t\right) =0 & \text{in }\left( 0,\infty \right), \vspace{0.2cm}\\
\varphi\left( 1,t\right) =\,\psi_{x} \left( 1,t\right) =\,w_{x}\left(
1,t\right) =0 & \text{in }\left( 0,\infty \right).
\end{array}
\right.  \label{cdt_1}
\end{equation}
The functions $\varphi,\,\psi$ and $w$ model, respectively, the vertical, shear angle and longitudinal displacements of the filament.
The coefficients $\rho_{1},\,\rho_{2} ,\, b,\,k,\,k_0 ,\,l$ and $\delta$ are positive constants. The unique dissipation considered in $(\ref{syst1})$ is played by the linear frictional damping $\delta \varphi _{t}$ (it is well known that, when $\delta=0$, $(\ref{syst1})$ is a conservative system).
\vskip0,1truecm
The Bresse system \cite{bres} has attracted the attention of many researchers for the last few years. Under different types of dissipative mechanisms, various stability results (in either bounded or unbounded domain) have been obtained. For this purpose, see \cite{afas}, \cite{alab}, \cite{alve}, \cite{char}, \cite{fato2}, \cite{gksh}, \cite{nowe}, \cite{sori2}, \cite{sori1}, \cite{souf} and \cite{wehb} in case of frictional dampings, and \cite{fato1}, \cite{kama}, \cite{liu2} and \cite{nawe} when the Bresse system is indirectly damped via the coupling with some other equations. In all these works, the considered dissipation is effective at least on the longitudinal displacement (third equation) or on the shear angle displacement (second equation). 
\vskip0,1truecm
In case of bounded domain, it was proved that the exponential stabiliy holds under some restrictions on the velocities of waves propagations, and the plynomial stability is valid in general with a decay rate depending on the regularity of initial data; see for example \cite{alab}, \cite{alve}, \cite{char}, \cite{fato2}, \cite{fato1}, \cite{kama}, \cite{liu2}, \cite{nawe}, 
\cite{nowe}, \cite{sori2}, \cite{sori1} and \cite{wehb}.
\vskip0,1truecm
However, when the domain is the whole line $\mathbb{R}$, the situation is completely different in the sens that no exponential stability result can be obtained and only polynomial stability results are proved on the $L^2$-norm of solutions (under some assumptions on the coefficients of the system) with a decay rate that can be improved by taking initial data more regular; see \cite{gksh} (one frictional damping acting on the third equation) and \cite{souf} (two frictional dampings acting on the second and third equations). When the frictional
damping is acting on the second equation, the authors of \cite{gksh} proved that there is no decay of solutions at all. 
\vskip0,1truecm
As far as we know, when only the veridical displacement (first equation) of Bresse system is damped via a frictional damping (the other two equations are totally free), the stability of Bresse system has never been considered in the literature. Unlike the papers cited above concerned with the case of bounded domain, we prove that, despite the presence of the linear frictional damping $\delta \varphi _{t}$, $(\ref{syst1})$ is never exponentially stable independently from the values of the coefficients of $(\ref{syst1})$. Moreover, we show that $(\ref{syst1})$ is at least polynomially stable with decay rates that can be improved by considering more smooth initial data.
\vskip0,1truecm
Our results show that the exponential stability of the overall Bresse system can not be guaranteed by a frictional dissipation working only in the veridical displacement. In comparaison with the known results cited above, this phenomenon means that Bresse system in a bounded domain is more dominated by its longitudinal and shear angle displacements than by its veridical displacement. 
\vskip0,1truecm
The paper is organized as follows: in section 2, we establish the existence, uniqueness and smoothness of solutions of $(\ref{syst1})-(\ref{cdt_1})$. In section 3, we show the non-exponential stability of $(\ref{syst1})-(\ref{cdt_1})$. Finally, the proof of the polynomial stability of $(\ref{syst1})-(\ref{cdt_1})$ is given in section 4.

\section{The semigroup setting}

In this section, we study the existence, uniqueness and smoothness of solutions for $(\ref{syst1})-\left( \ref{cdt_1}\right)$ using semigroup techniques. For this, let us considere the space
\begin{equation*}
\mathcal{H}=H_0^{1}\left( 0,1\right) \times L^{2}\left( 0,1\right)
\times H_*^{1}\left( 0,1\right) \times L_*^{2}\left(0,1\right) \times H_*^{1}\left( 0,1\right) \times L_*^{2}\left( 0,1\right),
\end{equation*}
where
\begin{equation*}
L_*^2 (0,1) =\left\{v\in L^2  (0,1),\,\int_0^1 v \,dx=0\right\}\quad\hbox{and}\quad H^{1}_* (0,1) =H^{1} (0,1)\cap L_*^2 (0,1) ,
\end{equation*}
equipped with the inner product
\begin{equation*}
\left\langle (\varphi_1 ,\,\tilde{\varphi}_1,\,\psi_1 ,\,\tilde{\psi}_1, \,w_1 ,\,\tilde{w}_1 )^T,(\varphi_2 ,\,\tilde{\varphi}_2,\,\psi_2 ,\,\tilde{\psi}_2, \,w_2 ,\,\tilde{w}_2 )^T\right\rangle_{\mathcal{H}}
\end{equation*} 
\begin{equation*}
\begin{array}{lll}
&=& k\left\langle \left( \varphi _{1x}+\psi_{1}+l\,w_{1}\right) ,\left(
\varphi _{2x}+\psi _{2}+l\,w_{2}\right) \right\rangle_{L^{2}\left(
0,1\right) }+k_{0}\left\langle \left( w_{1x}-l\varphi _{1}\right) ,\left(w_{2x}-l\varphi _{2}\right) \right\rangle _{L^{2}\left( 0,1\right)}\vspace{0.2cm}\\
&& \mbox{} +b\left\langle \psi _{1x},\psi _{2x}\right\rangle _{L^{2}\left( 0,1\right)} +\rho _{1}\left\langle \tilde{\varphi}_{1},\tilde{\varphi}_{2}\right\rangle _{L^{2}\left( 0,1\right)}+\rho _{2} \left\langle \tilde{\psi}_{1},\tilde{\psi}_{2}\right\rangle _{L^{2}\left( 0,1\right)}+\rho _{1}\left\langle \tilde{w}_{1},\tilde{w}_{2}\right\rangle _{L^{2}\left( 0,1\right) }.
\end{array}
\end{equation*}
Notice that, using the definition of $H^{1}_0 (0,1)$ and $H^{1}_* (0,1)$, we remark that, if 
\begin{equation*}
\left(\varphi ,\psi ,w\right)\in H^{1}_0 (0,1)\times H^{1}_* (0,1)\times H^{1}_* (0,1)
\end{equation*}
satisfying
\begin{equation*}
k\left\Vert \varphi_{x}+\psi +l\,w \right\Vert _{L^{2}\left( 0,1\right)}^{2} +b\left\Vert \psi _{x} \right\Vert_{L^{2}\left( 0,1\right) }^{2} +k_{0}\left\Vert w_{x}-l\varphi\right\Vert_{L^{2}\left( 0,1\right) }^{2}  =0,
\end{equation*}
then 
\begin{equation*}
\psi=0 ,\quad \varphi (x)=-c\sin\,(lx)\quad\hbox{and}\quad w(x)=c\cos\,(lx),
\end{equation*}
where $c$ is a constant such that $c=0$ or $l= m_0\pi$, for some $m_0\in \mathbb{Z}$. Furthermore, by assuming that 
\begin{equation}
l\ne m\pi,\quad\forall m\in \mathbb{Z} ,\label{pi}
\end{equation}
we get $\varphi = \psi =w=0$, and so, $\mathcal{H}$ is a Hilbert space with respect to the generated norm
\begin{equation*}
\begin{array}{lll}
\left\Vert \left( \varphi ,\,\tilde{\varphi},\,\psi ,\,\tilde{\psi},\,w,\,\tilde{w}\right)^{T} \right\Vert _{\mathcal{H}}^2 &=& k\left\Vert \varphi_{x}+\psi +l\,w \right\Vert _{L^{2}\left( 0,1\right)}^{2} +b\left\Vert \psi _{x}\right\Vert_{L^{2}\left( 0,1\right) }^{2} +k_{0}\left\Vert w_{x}-l\varphi \right\Vert_{L^{2}\left( 0,1\right) }^{2}\vspace{0.2cm}\\
&& \mbox{} +\rho _{1}\left\Vert \tilde{\varphi}\right\Vert_{L^{2}\left( 0,1\right) }^{2}+\rho _{2}\Vert \tilde{\psi}\Vert_{L^{2}\left( 0,1\right) }^{2} +\rho _{1}\left\Vert \tilde{w}\right\Vert_{L^{2}\left( 0,1\right)}^{2} .
\end{array}
\end{equation*}
The definition of $L_*^2 (0,1)$ allows to apply Poincar\'e's inequality in 
$H_*^1 (0,1)$, and the property 
\begin{equation*}
\int_0^1 v(x)dx =0,
\end{equation*}
for $v\in\{\psi ,w\}$ can be assumed without lose of generality thanks to a classical change of variables; see, for example, Remark 2.1 of \cite{guki}.
\vskip0,1truecm
Now, we consider the vectors
\begin{equation*}
\Phi =\left( \varphi ,\,\tilde{\varphi},\,\psi ,\,\tilde{\psi},\,w,\,\tilde{w}\right)^{T} \,\,\hbox{and}\,\,\Phi_0 =
\left( \varphi_{0},\,\varphi_{1},\,\psi_{0},\,\psi_{1},\,w_{0},\,w_{1}\right)^{T} ,
\end{equation*}
where $\tilde{\varphi}=\varphi_{t}$, $\tilde{\psi}=\psi_{t}$ and $\tilde{w}=w_{t}$. System $(\ref{syst1})-\left( \ref{cdt_1}\right)$ can be formulated as the following first order system: 
\begin{equation}
\left\{
\begin{array}{ll}
\Phi_{t}=\mathcal{A}\Phi\quad\quad\quad\text{in } \left( 0,\infty \right) ,\vspace{0.2cm}\\
\Phi \left( 0\right) =\Phi_{0} ,
\end{array}
\right.  \label{syst_2}
\end{equation}
where 
\begin{equation}
\mathcal{A}\Phi =\left(
\begin{array}{c}
\tilde{\varphi} \vspace{0.2cm}\\
\dfrac{k}{\rho _{1}}\left( \varphi _{x}+\psi +l\,w\right) _{x}+\dfrac{lk_{0}%
}{\rho _{1}}\left( w_{x}-l\varphi \right) -\dfrac{\delta }{\rho _{1}}\tilde{\varphi} \vspace{0.2cm}\\
\tilde{\psi} \vspace{0.2cm}\\
\dfrac{b}{\rho _{2}}\psi _{xx}-\dfrac{k}{\rho _{2}}\left( \varphi _{x}+\psi
+l\,w\right) \vspace{0.2cm}\\
\tilde{w} \vspace{0.2cm}\\
\dfrac{k_{0}}{\rho _{1}}\left( w_{x}-l\varphi \right) _{x}-\dfrac{lk}{\rho
_{1}}\left( \varphi _{x}+\psi +l\,w\right) 
\end{array}
\right)  \label{oper_A}
\end{equation}
with domain 
\begin{equation*}
D\left( \mathcal{A}\right) =\left\{
\begin{array}{c}
\Phi \in \mathcal{H}\mid \,\varphi \in \,H^{2}\left( 0,1\right)\cap H_0^{1}\left( 0,1\right) 
;\,\psi ,\,w\in \,H^{2}\left( 0,1\right)\cap H_*^{1}\left( 0,1\right) ; \\
\tilde{\varphi}\in H_0^{1}\left( 0,1\right);\,\tilde{\psi},\,\tilde{w}\in \,H_*^{1}\left( 0,1\right) ;\,\,\psi_{x}\left( 0\right) =w_{x}\left( 0\right) =\psi_{x}\left( 1\right) =w_{x}\left( 1\right) =0
\end{array}
\right\} .
\end{equation*}
\vskip0,2truecm
\begin{theorem}\label{Theorem 1.1}
Assume that $(\ref{pi})$ holds. Then, for any $m\in \mathbb{N}$ and $\Phi_0 \in D(\mathcal{A}^m)$, system $(\ref{syst_2})$ admits a unique solution
\begin{equation}
\Phi\in \cap_{j=0}^m C^{m-j} \left(\mathbb{R}_{+} ;D\left(\mathcal{A}^j\right)\right). \label{exist}
\end{equation}
\end{theorem}
\vskip0,2truecm
\begin{proof} We remark that $D(\mathcal{A})$ is dense in $\mathcal{H}$. Now, direct calculation gives
\begin{equation}
\left\langle \mathcal{A}\Phi ,\Phi \right\rangle _{\mathcal{H}}=-\delta
\left\Vert \tilde{\varphi}\right\Vert _{L^{2}\left( 0,1\right) }^{2}\leq 0. \label{dissp}
\end{equation}
Hence, $\mathcal{A}$ is a dissipative operator. 
\vskip0,1truecm
Next, we show that $0\in\rho\left( \mathcal{A}\right)$. Let $F=(f_1 ,\cdots,f_6)^T \in \mathcal{H}$. We prove that there exists $Z=(z_1 ,\cdots,z_6)^T \in D\left( \mathcal{A}\right)$ satisfying
\begin{equation}
\mathcal{A} Z=F.  \label{ZF}
\end{equation}
Indeed, first, the first, third and fifth equations in $(\ref{ZF})$ are equivalent to
\begin{equation}
z_2 =f_1 ,\quad z_4 =f_3 \quad\hbox{and}\quad z_6 =f_5 ,  \label{z1f1}
\end{equation}
and then
\begin{equation}
z_2 \in H_0^{1} \left( 0,1\right) \quad\hbox{and}\quad z_4 ,\,z_6 \in H_{*}^{1}\left( 0,1\right). \label{z2f2}
\end{equation}
Second, substitute $z_2$ into the second equation in $(\ref{ZF})$, we find that the second, fourth and sixth equations in $(\ref{ZF})$
are reduced to
\begin{equation}
\left\{
\begin{array}{ll}
k\left( z_{1x}+z_{3} +l\,z_{5}\right)_{x} +lk_{0}\left( z_{5x}-lz_{1} \right) =\delta f_{1} +\rho_1 f_2 , \vspace{0.2cm}\\
bz_{3xx} -k\left( z_{1x}+z_{3} +l\,z_{5}\right)=\rho_2 f_4 , \vspace{0.2cm}\\
k_{0}\left( z_{5x}-lz_{1} \right)_{x} -lk\left( z_{1x}+z_{3} +l\,z_{5}\right) =\rho_1 f_6.
\end{array}
\right. \label{z5f5}
\end{equation}
To prove that $(\ref{z5f5})$ admits a solution $(z_1 ,z_3 ,z_5 )$ satisfying
\begin{equation}
\left\{
\begin{array}{ll}
z_1 \in H^{2} \left( 0,1\right)\cap H_0^{1} \left( 0,1\right),\quad z_3 ,\,z_5 \in H^{2} \left( 0,1\right)\cap H_*^{1}\left( 0,1\right) , \vspace{0.2cm}\\
z_{3x} (0)=z_{5x} (0)=z_{3x} (1)=z_{5x} (1)=0, \label{z6f6}
\end{array}
\right. 
\end{equation}
we define the following bilinear form:
\begin{equation*}
a \left(\left(v_1 ,\,v_2 ,\,v_3 \right) , \left(w_1 ,\,w_2 ,\,w_3\right) \right) =k\left\langle v_{1x} +v_{2}+lv_{3} ,w_{1x}+w_{2}+lw_{3} \right\rangle_{L^{2}\left(0,1\right)}
\end{equation*}
\begin{equation*}
+b\left\langle v_{2x},w_{2x}\right\rangle_{L^{2}\left( 0,1\right)} +k_{0}\left\langle v_{3x}-lv_{1}, w_{3x}-lw_{1}\right\rangle_{L^{2}\left( 0,1\right)} ,\quad\,\,\forall \left(v_1 ,\,v_2 ,\,v_3\right)^T ,\,\left(w_1 ,\,w_2 ,\,w_3\right)^T\in \mathcal{H}_0 ,
\end{equation*}
and the following linear form:
\begin{equation*}
L \left( v_1 ,\,v_2 ,\,v_3\right) =\left\langle \delta f_{1} +\rho_1 f_2 ,v_{1} \right\rangle_{L^{2}\left(0,1\right)}+\left\langle \rho_2 f_4 ,v_{2} \right\rangle_{L^{2}\left(0,1\right)}
+\left\langle \rho_1 f_6 , v_{3} \right\rangle_{L^{2}\left(0,1\right)},\quad\forall \left(v_1 ,\,v_2 ,\,v_3\right)^T \in\mathcal{H}_0 ,
\end{equation*}
where
\begin{equation*}
\mathcal{H}_0 =H_{0}^{1} \left( 0,1\right) \times H_{*}^{1}\left( 0,1\right)\times H_{*}^{1}\left( 0,1\right) .
\end{equation*}
Thus, the variational formulation of $(\ref{z5f5})$ is given by
\begin{equation}
a \left(\left(z_1 ,\,z_3 ,\,z_5 \right), \left(w_1 ,\,w_2 ,\,w_3\right)\right)=L \left(w_1 ,\,w_2 ,\,w_3\right),\,\forall \left(w_1 ,\,w_2 ,\,w_3\right)^T \in \mathcal{H}_0 . \label{z7f7}
\end{equation}
From the Lax-Milgram theorem, it follows that $(\ref{z7f7})$ has a unique solution
\begin{equation*}
\left(z_1 ,\,z_3 ,\,z_5 \right)\in\mathcal{H}_0 .
\end{equation*}
Therefore, using classical elliptic regularity arguments, we conclude that $\left(z_1 ,\,z_3 ,\,z_5 \right)$ solves $(\ref{z5f5})$ and satisfies the regularity and boundary conditions $(\ref{z6f6})$. This proves that $(\ref{ZF})$ has a unique solution $Z\in D\left( \mathcal{A}\right)$. By the resolvent identity, we have $\lambda I -\mathcal{A}$ is surjective, for any $\lambda >0$ (see \cite{liu1}), where $I$ denotes the identity operator. Consequently, the Lumer-Phillips theorem implies that $\mathcal{A}$ is the infinitesimal generator of a linear $C_{0}$ semigroup of contractions on $\mathcal{H}$. So, Theorem \ref{Theorem 1.1} holds (see \cite{pazy}).
\end{proof} 
\vskip0,1truecm
The proof of the non-exponential and polynomial stability for $(\ref{syst_2})$ is based on the following two frequency domain theorems:
\vskip0,2truecm
\begin{theorem}\label{Theorem 2.1}
(\cite{huan} and \cite{prus}) A $C_{0}$ semigroup of contractions on a Hilbert space $\mathcal{H}$ generated by an operator
$\mathcal{A}$ is exponentially stable if and only if
\begin{equation}
i\mathbb{R} \subset \rho\left( \mathcal{A}\right)\quad\hbox{and}\quad \sup_{\lambda\in \mathbb{R}}\left\Vert \left( i\lambda I-\mathcal{A}\right) ^{-1} \right\Vert_{\mathcal{L}\left( \mathcal{H}\right) } <\infty. \label{expon}
\end{equation}
\end{theorem}
\vskip0,2truecm
\begin{theorem}\label{Theorem 2.2}
(\cite{liu0}) If a bounded $C_{0}$ semigroup $e^{t\mathcal{A}}$ on a Hilbert space $\mathcal{H}$ generated by an operator
$\mathcal{A}$ satisfies, for some $j\in \mathbb{N}^*$,
\begin{equation}
i\mathbb{R} \subset \rho\left( \mathcal{A}\right)\quad\hbox{and}\quad \sup_{\vert\lambda\vert\geq 1}\dfrac{1}{\lambda^j}\left\Vert \left( i\lambda I-\mathcal{A}\right) ^{-1} \right\Vert_{\mathcal{L}\left( \mathcal{H}\right) } <\infty. \label{polyc}
\end{equation}
Then, for any $m\in \mathbb{N}^*$, there exists a positive constant $c_m$ such that
\begin{equation}
\left\Vert e^{t\mathcal{A}} z_0\right\Vert_{\mathcal{H}} \leq c_m  \left\Vert z_0 \right\Vert_{D\left(\mathcal{A}^m\right)} \left(\dfrac{\ln\,t}{t}\right)^{\dfrac{m}{j}}\ln\,t,\quad\forall z_0\in D\left(\mathcal{A}^m\right),\,\,\forall t>0. \label{polys}
\end{equation}
\end{theorem}

\section{Lack of exponential stability of $(\ref{syst1})-(\ref{cdt_1})$}

In this section, we prove our first main result which is stated as follows:
\vskip0,2truecm
\begin{theorem}\label{Theorem 3.1}
Assume that $(\ref{pi})$ holds. Then the semigroup associated with $\left( \ref{syst_2}\right)$ is not exponentially stable.
\end{theorem}
\vskip0,2truecm
\begin{proof} Using Theorem \ref{Theorem 2.1}, it is enough to prove that the second condition in $\left( \ref{expon}\right)$ is not satisfied. To do so, we prove that there exists a sequence $(\lambda_{n})_n \subset \mathbb{R}$ such that
\begin{equation*}
\lim_{n\to \infty}\left\Vert \left( i\lambda _{n}I-\mathcal{A}\right) ^{-1}\right\Vert_{\mathcal{L}\left( \mathcal{H}\right) } =\infty,
\end{equation*}
which is equivalent to find a sequence $(F_{n})_n\subset \mathcal{H}$ satisfying
\begin{equation}
\left\Vert F_{n}\right\Vert_{\mathcal{H}}\leq 1,\quad\forall n\in\mathbb{N} \label{Fn}
\end{equation}
and
\begin{equation}
\lim_{n\to \infty}\Vert\left( i\lambda _{n}I-\mathcal{A}\right)^{-1}F_{n}\Vert_{\mathcal{H}} = \infty.  \label{eq_3}
\end{equation}
For this purpose, let
\begin{equation*}
\Phi_{n}=\left( i\lambda _{n}I-\mathcal{A}\right)^{-1} F_n , \quad\forall n\in\mathbb{N} ,
\end{equation*}
where
\begin{equation*}
\Phi _{n} =\left(\varphi_n ,{\tilde{\varphi}}_n ,\psi_n ,{\tilde{\psi}}_n ,w_n ,{\tilde{w}}_n \right)^T\quad\hbox{and}\quad F_n =\left(f_{1n} ,\cdots,f_{6n} \right)^T.
\end{equation*}
Then we have to find $(\lambda_{n})_n\subset \mathbb{R}$, $(F_{n})_n\subset \mathcal{H}$ and $(\Phi_{n})_n\subset D(\mathcal{A})$ satisfying 
$\left( \ref{Fn}\right)$,
\begin{equation}
\lim_{n\to \infty}\Vert\Phi _{n}\Vert_{\mathcal{H}}=\infty\quad\hbox{and}\quad i\lambda _{n}\Phi _{n}-\mathcal{A}\Phi _{n}=F_{n} ,  \,\,\forall n\in\mathbb{N}. \label{eq_4}
\end{equation}
The equation in $\left( \ref{eq_4}\right)$ is equivalent to  
\begin{equation}
\left\{
\begin{array}{l}
i\lambda _{n}\varphi _{n}-{\tilde{\varphi}}_{n}=f_{1n} ,\vspace{0.2cm}\\
i\rho _{1}\lambda _{n}{\tilde{\varphi}}_{n}-k\left( \varphi _{nx}+\psi _{n}+l\,w_{n}\right)
_{x}-lk_{0}\left( w_{nx}-l\varphi _{n}\right) +\delta {\tilde{\varphi}}_{n} =\rho_1 f_{2n} ,\vspace{0.2cm}\\
i\lambda _{n}\psi _{n}-{\tilde{\psi}}_{n}\;=f_{3n} ,\vspace{0.2cm}\\
i\rho _{2}\lambda _{n}{\tilde{\psi}}_{n}-b\psi _{nxx}+k\left( \varphi _{nx}+\psi
_{n}+l\,w_{n}\right) =\rho_2 f_{4n} ,\vspace{0.2cm}\\
i\lambda _{n}w_{n}-{\tilde{w}}_{n}=f_{5n} ,\vspace{0.2cm}\\
i\rho _{1}\lambda _{n}{\tilde{w}}_{n}-k_{0}\left( w_{nx}-l\varphi _{n}\right)
_{x}+lk\left( \varphi _{nx}+\psi _{n}+l\,w_{n}\right) =\rho_1 f_{6n} .
\end{array}
\right.  \label{eq_2_5}
\end{equation}
Choosing
\begin{equation}
f_{1n}=f_{3n}=f_{5n} =0. \label{fn0}
\end{equation}
Then system $\left( \ref{eq_2_5}\right)$ becomes
\begin{equation}
\left\{
\begin{array}{l}
{\tilde{\varphi}}_{n} =i\lambda_{n}\varphi_{n} ,\quad {\tilde{\psi}}_{n} =i\lambda_{n}\psi_{n} ,\quad {\tilde{w}}_{n} =i\lambda _{n} w_n ,\vspace{0.2cm}\\
\left(i\delta \lambda_{n}-\rho _{1}\lambda _{n}^{2}\right)\varphi _{n}-k\left( \varphi _{nx}+\psi_{n}+l\,w_{n}\right) _{x}-lk_{0}\left( w_{nx}-l\varphi _{n}\right) =\rho_1 f_{2n} ,\vspace{0.2cm}\\
-\rho _{2}\lambda _{n}^{2}\psi _{n}-b\psi _{nxx}+k\left( \varphi
_{nx}+\psi _{n}+l\,w_{n}\right) =\rho_2 f_{4n} ,\vspace{0.2cm}\\
-\rho _{1}\lambda _{n}^{2}w_{n}-k_{0}\left( w_{nx}-l\varphi _{n}\right)
_{x}+lk\left( \varphi _{nx}+\psi _{n}+l\,w_{n}\right) =\rho_1 f_{6n} . 
\end{array}
\right.  \label{eq_2_6}
\end{equation}
To simplify the calculations, we put $N=n\pi$. Some of the computations below were done in \cite{afas}. Now, we consider three cases.
\vskip0,1truecm
{\bf Case 1}: $\dfrac{b}{\rho _{2}}=\dfrac{k_{0}}{\rho _{1}}$. We choose
\begin{equation}
\left\{
\begin{array}{l}
\varphi_{n} ={\tilde{\varphi}}_{n} =0, \vspace{0.2cm}\\
\psi_n (x)=\alpha_1\cos\, \left(Nx\right),\quad {\tilde{\psi}}_{n} (x)=i\alpha_1\lambda_n\cos\, \left(Nx\right), \vspace{0.2cm}\\ 
w_n (x)= \alpha_2\cos\, \left(Nx\right),\quad {\tilde{w}}_{n} (x)=i\alpha_2\lambda_n \cos\, \left(Nx\right) , 
\end{array}
\right. \label{fn2}
\end{equation}
\begin{equation}
f_{2n}=0, \quad f_{4n} (x)=-\frac{lk_{0}}{\rho _{2}}\alpha_2\cos\,\left(Nx\right),\quad f_{6n} (x)=-\frac{l^{2}k_{0}}{\rho _{1}}\alpha_2\cos\, \left(Nx\right) \label{fn1}
\end{equation}
and 
\begin{equation}
\lambda _{n} =N\sqrt{\dfrac{k_{0}}{\rho _{1}}}, \label{fn7}
\end{equation}
where $\alpha_1,\,\alpha_2\in \mathbb{R}$. We have $\Phi_{n} \in D(\mathcal{A})$ and $F_{n} \in \mathcal{H}$. On the other hand,
$\left( \ref{eq_2_6}\right)$ is satisfied if and only if
\begin{equation}
\left\{
\begin{array}{l}
k\alpha_1 +l\left( k+k_{0}\right) \alpha_2 =0, \vspace{0.2cm}\\
\left[ -\lambda _{n}^{2}+\dfrac{b}{\rho _{2}} N^{2}+\dfrac{k}{\rho _{2}}\right] \alpha_1 +\dfrac{lk}{\rho _{2}}\alpha_2 =-\dfrac{lk_{0}}{\rho _{2}} \alpha_2, \vspace{0.2cm}\\
\dfrac{lk}{\rho _{1}} \alpha_1 +\left[ -\lambda _{n}^{2}+\dfrac{k_{0}}{\rho _{1}} N^{2}+\dfrac{l^{2}k}{\rho_{1}}\right] \alpha_2 =-\dfrac{l^{2}k_{0}}{\rho _{1}} \alpha_2.
\end{array}
\right. \label{fn4}
\end{equation}
According to $\left( \ref{fn7}\right)$ and because $\dfrac{b}{\rho _{2}}=\dfrac{k_{0}}{\rho _{1}}$, we have 
\begin{equation*}
-\lambda _{n}^{2}+\dfrac{k_0}{\rho _{1}} N^{2}=-\lambda _{n}^{2}+\dfrac{b}{\rho _{2}} N^{2}=0,
\end{equation*}%
and therefore, the system $\left(\ref{fn4}\right)$ is equivalent to
\begin{equation}
\alpha_1 =-l\left( 1+\frac{k_{0}}{k}\right) \alpha_2 . \label{fn5}
\end{equation}
Choosing \begin{equation*}
\alpha_2 =\frac{\rho _{1} \rho _{2}}{lk_{0} \sqrt{\rho _{1}^2 +l^2\rho _{2}^{2} }}
\end{equation*}
and using $\left( \ref{fn0}\right)$ and $\left( \ref{fn1}\right)$, we obtain
\begin{equation*}
\left\Vert F_{n}\right\Vert _{\mathcal{H}}^{2} =\left\Vert f_{4n}\right\Vert _{L^{2} \left( 0,1\right)}^{2}+\left\Vert f_{6n}\right\Vert_{L^{2}\left( 0,1\right)}^{2} =\left( \frac{lk_{0}}{\rho _{2}}\right) ^{2}\left[ 1+\left(\frac{l\rho _{2}}{\rho _{1}}\right) ^{2}\right] \alpha_2^{2}\int_{0}^{1}\cos ^{2} \left(Nx\right)\,dx
\end{equation*}
\begin{equation*}
\leq\left( \frac{lk_{0}}{\rho _{2}}\right) ^{2}\left[ 1+\left( \frac{l\rho _{2}}{\rho _{1}}\right) ^{2}\right] \alpha_2^{2} = 1 .
\end{equation*}%
On the other hand, we have  
\begin{equation*}
\left\Vert \Phi _{n}\right\Vert _{\mathcal{H}}^{2} \geq k_0\left\Vert w_{nx} -l\varphi_n\right\Vert _{L^{2}\left( 0,1\right) }^{2}
= k_0\left\Vert w_{nx}\right\Vert _{L^{2}\left( 0,1\right) }^{2} =\frac{k_0\alpha_2^{2}}{2} N^{2}\int_{0}^{1} \left[1-\cos\,\left(2Nx\right)\right]\,dx=\frac{k_0\alpha_2^{2}}{2} N^{2} ,
\end{equation*}
hence
\begin{equation}
\lim_{n\to\infty}\left\Vert \Phi _{n}\right\Vert _{\mathcal{H}} = \infty. \label{fn69}
\end{equation}
\vskip0,1truecm
{\bf Case 2}: $\dfrac{b}{\rho _{2}}\ne\dfrac{k_{0}}{\rho _{1}}$ and $k\ne k_0$. We choose
\begin{equation}
f_{2n} = f_{4n} =0,\quad f_{6n} (x) =\cos\, \left(Nx\right), \label{fn025}
\end{equation}
\begin{equation}
\left\{
\begin{array}{l}
\varphi _{n} (x) =\alpha_1 \sin\, \left(Nx\right), \quad {\tilde{\varphi}}_n (x)=i\alpha_1\lambda_n \sin\, \left(Nx\right),\vspace{0.2cm}\\
\psi _{n} (x)=\alpha_2\cos\, \left(Nx\right),\quad {\tilde{\psi}}_n (x)=i\alpha_2\lambda_n \cos\, \left(Nx\right),\vspace{0.2cm}\\
w_{n} (x)=\alpha_3\cos\, \left(Nx\right),\quad {\tilde{w}}_n (x)=i\alpha_3\lambda_n\cos\, \left(Nx\right)
\end{array}
\right. \label{fn015}
\end{equation}
and
\begin{equation}
\lambda_n =\sqrt{\dfrac{k_0}{\rho_1} N^2 +\dfrac{l^2 k}{\rho_1}} ,\label{fn0650+}
\end{equation}
where $\alpha_1,\,\alpha_2 ,\,\alpha_3\in \mathbb{C}$ depending on $N$. Notice that, according to these choices, $\Phi_{n} \in D(\mathcal{A})$, $F_{n} \in \mathcal{H}$ and 
\begin{equation}
\left\Vert F_{n}\right\Vert _{\mathcal{H}}^{2} =\left\Vert f_{6n}\right\Vert_{L^{2} \left( 0,1\right)}^{2} =\int_{0}^{1}\cos ^{2} \left(Nx\right)\,dx\leq 1 . \label{fn065+}
\end{equation} 
On the other hand, thanks to $\left( \ref{fn0}\right)$, $\left( \ref{fn025}\right)$ and $\left( \ref{fn015}\right)$, the first three equations in $\left( \ref{eq_2_6}\right)$ are satisfied, and the last three ones are equivalent to
\begin{equation}
\left\{
\begin{array}{l}
\left[\left(k-\mu_n\right)N^2 -\rho _{1}\lambda _{n}^{2}+l^2 k_0 \right]\alpha_1 +kN\alpha_2 +l\left(k+k_0\right)N\alpha_3 =0, \vspace{0.2cm}\\
kN\alpha_1 +\left(bN^2 -\rho _{2}\lambda _{n}^{2}+k\right)\alpha_2 +kl\alpha_3 =0,\vspace{0.2cm}\\
l\left(k+k_0\right)N\alpha_1 +lk\alpha_2 +\left(k_0 N^2 -\rho _{1}\lambda _{n}^{2}+l^2 k\right)\alpha_3 =\rho_1 ,
\end{array}
\right.  \label{eq_2_61}
\end{equation}
where we note
\begin{equation}
\mu_n = \frac{-i\delta \lambda_n}{N^2} . \label{mun*}
\end{equation}
From the choice $\left( \ref{fn0650+}\right)$, we see that the last equation in $\left( \ref{eq_2_61}\right)$ is equivalent to 
\begin{equation}
\alpha_2 =-\dfrac{k+k_0}{k} N\alpha_1 +\dfrac{\rho_1}{lk},\label{eq_2_63}
\end{equation}
so, substituting in the first two equations in $\left( \ref{eq_2_61}\right)$, we get
\begin{equation}
\alpha_3 = a_1 N\alpha_{1} +a_2\label{eq_2_66}
\end{equation}
and
\begin{equation}
\alpha_1 = \frac{\left[l\left(k+k_0\right)a_2 +\dfrac{\rho_1}{l}\right]N}{\left[2k_0 +\mu_n -l\left(k+k_0\right)a_1 \right]N^2 +l^2 \left(k-k_0\right)},\label{eq_2_67}
\end{equation}
where
\begin{equation*}
\left\{
\begin{array}{l}
a_1 =\frac{k+k_0}{lk^2}\left(b-\dfrac{\rho_2 k_0}{\rho_1}\right)N^2 +\dfrac{k_0}{lk}-\dfrac{l\rho_2 (k+k_0)}{\rho_1 k} , \vspace{0.2cm}\\
a_2 =\frac{\rho_1}{(lk)^2}\left[\left(\dfrac{\rho_2 k_0}{\rho_1} -b\right)N^2 +\dfrac{l^2 \rho_2 k}{\rho_1}-k\right] .
\end{array}
\right. 
\end{equation*} 
To simplify the computations, we put
\begin{equation*}
\left\{
\begin{array}{l}
a_3 =\dfrac{\rho_1\left(k+k_0\right)}{lk^2}\left(\dfrac{\rho_2 k_0}{\rho_1} -b\right) , \quad
a_4 =\dfrac{\left(k+k_0\right)^2}{k^2}\left(\dfrac{\rho_2 k_0}{\rho_1} -b\right) , \vspace{0.2cm}\\
a_5 = \frac{l\rho_2\left(k+k_0\right) }{k} -\dfrac{k_0 \rho_1}{lk}, \quad
a_6 =\dfrac{l^2 \rho_2\left(k+k_0\right)^2}{\rho_1 k} +\dfrac{k_0 \left(k-k_0\right)}{k} 
\end{array}
\right. 
\end{equation*}
and
\begin{equation*}
\left\{
\begin{array}{l}
d_0 = \frac{k+k_0}{lk^2}\left(b-\dfrac{\rho_2 k_0}{\rho_1}\right), \quad
d_1 = \frac{\rho_1}{(lk)^2} \left(\dfrac{\rho_2 k_0}{\rho_1} -b\right), \vspace{0.2cm}\\
d_2 = \dfrac{k_0}{lk}-\dfrac{l\rho_2 (k+k_0)}{\rho_1 k} , \quad
d_3 = \frac{\rho_1}{l^2 k} \left(\dfrac{l^2 \rho_2 }{\rho_1}-1\right).
\end{array}
\right. 
\end{equation*}
Then
\begin{equation*}
N\alpha_1 = \frac{a_3 N^4 +a_5 N^2}{a_4 N^4 +\left(\mu_n +a_6\right)N^2 +l^2 \left(k-k_0\right)}
\end{equation*}
and (notice that $d_0 a_3 +d_1 a_4 =0$)
\begin{equation}
\alpha_3 = \frac{\left(d_0 N^2 +d_2\right)\left(a_3 N^4 +a_5 N^2\right)}{a_4 N^4 +\left(\mu_n +a_6\right)N^2 +l^2 \left(k-k_0\right)}+d_1 N^2 +d_3\label{eq_2_67+}
\end{equation}
\begin{equation*}
=\frac{\left(d_0 a_5 +d_2 a_3 +d_3 a_4 +d_1 a_6 +d_1 \mu_n\right)N^4 +\left(d_2 a_5 +d_3 a_6 + l^2 \left(k-k_0\right)d_1 +d_3 \mu_n\right) N^2 +l^2 \left(k-k_0\right)d_3 }{a_4 N^4 +\left(\mu_n +a_6\right)N^2 +l^2 \left(k-k_0\right)},
\end{equation*}
Because $\dfrac{b}{\rho _{2}}\ne\frac{k_{0}}{\rho _{1}}$ and $k\ne k_0$, then $a_4 \ne 0$ and
\begin{equation}
d_0 a_5 +d_2 a_3 +d_3 a_4 +d_1 a_6 = \frac{\rho_1}{(lk)^2} \left(\dfrac{\rho_2 k_0}{\rho_1} -b\right)\left(k_0 -k\right)\ne 0.\label{eq_2_67++}
\end{equation}
On the other hand, 
\begin{equation}
\lim_{n\to \infty}  \mu_n = \lim_{n\to \infty} \frac{-i\delta\lambda_n}{N^2} =\lim_{n\to \infty} \frac{-i\delta}{N^2}\sqrt{\dfrac{k_0}{\rho_1} N^2 +\dfrac{l^2 k}{\rho_1}} =0.  \label{mun} 
\end{equation}
Then we deduce from $\left(\ref{eq_2_67+}\right)$, $\left(\ref{eq_2_67++}\right)$ and $\left(\ref{mun}\right)$ that
\begin{equation}
\lim_{n\to \infty} \alpha_3 =\dfrac{d_0 a_5 +d_2 a_3 +d_3 a_4 +d_1 a_6}{a_4} \ne 0,  \label{eq_2_71} 
\end{equation}
hence
\begin{equation}
\lim_{n\to \infty} \vert\alpha_3\vert\lambda_n = \infty. \label{eq_2_71++} 
\end{equation}
Now, we have 
\begin{equation*}
\left\Vert \Phi _{n}\right\Vert _{\mathcal{H}}^{2} \geq\rho_1 \left\Vert {\tilde w}_{n}\right\Vert _{L^{2}\left( 0,1\right) }^{2} =\rho_1\left(\vert\alpha_3\vert \lambda_n\right)^{2} \int_{0}^{1} \cos^2\,\left(Nx\right)\,dx 
= \frac{\rho_1}{2}\left(\vert\alpha_3\vert \lambda_n\right)^{2} \int_{0}^{1} \left[1+\cos\,\left(2Nx\right)\right]\,dx=\frac{\rho_1}{2}\left(\vert\alpha_3\vert \lambda_n\right)^{2},
\end{equation*}
then by $\left( \ref{eq_2_71++}\right)$ we get $\left( \ref{fn69}\right)$.
\vskip0,1truecm
{\bf Case 3}: $\dfrac{b}{\rho _{2}}\ne\dfrac{k_{0}}{\rho _{1}}$ and $k= k_0$. We consider the choices $\left( \ref{fn0}\right)$,
\begin{equation}
\lambda_n =\sqrt{\dfrac{b}{\rho_2} N^2 +\dfrac{k}{2\rho_2}}, \label{fn0650++}
\end{equation}
\begin{equation}
f_{2n} = 0,\quad f_{4n} (x)= \alpha_2 C_n \cos\, \left(Nx\right) ,\quad f_{6n} (x) =\alpha_2 D_n\cos\, \left(Nx\right) \label{fn025++}
\end{equation}
and $\left( \ref{fn015}\right)$ with 
\begin{equation}
\alpha_1 = \left(\frac{\rho_1 D_n}{2lk} -\frac{1}{2}\right)\frac{\alpha_2}{N} \quad\hbox{and}\quad \alpha_3 =0, \label{fn025+++}
\end{equation}
where
\begin{equation*}
C_n =\frac{\rho_1}{2l\rho_2}D_n\quad\hbox{and}\quad D_n =\frac{2lk}{\rho_1}\left(\frac{1}{2} -\frac{k}{k+\frac{l^2 k}{N^2}-\mu_n -\frac{\rho_1 \lambda_n^2 }{N^2}}\right) . 
\end{equation*}
According to $\left( \ref{mun*}\right)$ and $\left( \ref{fn0650++}\right)$, we remark that $\lim_{n\to \infty} \mu_n =0$, and then 
\begin{equation*}
\lim_{n\to \infty} D_n = \frac{2lk}{\rho_1}\left(\frac{1}{2} -\frac{k}{k -\frac{\rho_1 b}{\rho_2}}\right)\quad\hbox{and}\quad 
\lim_{n\to \infty} C_n = \frac{k}{\rho_2}\left(\frac{1}{2} -\frac{k}{k -\frac{\rho_1 b}{\rho_2}}\right)
\end{equation*}
(these limits exist since $\dfrac{b}{\rho _{2}}\ne\dfrac{k_{0}}{\rho _{1}}$ and $k=k_0$), so, the sequence $\left(\vert C_n\vert^2 +\vert D_n\vert^2\right)_n$ is bounded. Then we choose 
\begin{equation}
\alpha_2 =\frac{1}{\sqrt{sup_{n\in \mathbb{N}}\left(\vert C_n\vert^2 +\vert D_n\vert^2\right)}} . \label{alpha2} 
\end{equation}
According to these choices, we see that $\Phi_{n} \in D(\mathcal{A})$, $F_{n} \in \mathcal{H}$ and, using $\left( \ref{fn0}\right)$, $\left( \ref{fn025++}\right)$ and $\left( \ref{alpha2}\right)$, we find  
\begin{equation*}
\left\Vert F_{n}\right\Vert _{\mathcal{H}}^{2} =\left\Vert f_{4n}\right\Vert_{L^{2} \left( 0,1\right)}^{2}+\left\Vert f_{6n}\right\Vert_{L^{2} \left( 0,1\right)}^{2} = \left(\vert C_n\vert^2 +\vert D_n\vert^2\right)\alpha_2^2
\int_{0}^{1}\cos ^{2} \left(Nx\right)\,dx 
\leq \left(\vert C_n\vert^2 +\vert D_n\vert^2\right)\alpha_2^2 \leq 1.
\end{equation*}
On the other hand, thanks to $\left( \ref{fn0}\right)$, $\left( \ref{fn015}\right)$ and $\left( \ref{fn025++}\right)$, the first three equations in $\left( \ref{eq_2_6}\right)$ are satisfied, and because $k=k_0$ and $\alpha_3 =0$, the last three equations in $\left( \ref{eq_2_6}\right)$ are equivalent to
\begin{equation}
\left\{
\begin{array}{l}
\left[\left(k-\mu_n\right)N^2 -\rho _{1}\lambda _{n}^{2}+l^2 k \right]\alpha_1 +kN\alpha_2 =0, \vspace{0.2cm}\\
kN\alpha_1 +\left(bN^2 -\rho _{2}\lambda _{n}^{2}+k\right)\alpha_2 = \rho_2\alpha_2 C_n,\vspace{0.2cm}\\
2lkN\alpha_1 +lk\alpha_2 =\rho_1 \alpha_2 D_n.
\end{array}
\right.  \label{eq_2_61+++}
\end{equation}
The first equation in $\left( \ref{eq_2_61+++}\right)$ is satisfied thanks to the definition of $\alpha_1$ and $D_n$, the second equation in $\left( \ref{eq_2_61+++}\right)$ holds according to the definition of $\lambda_n$, $\alpha_1$ and $C_n$, and the last equation in $\left( \ref{eq_2_61+++}\right)$ is satisfied from the definition of $\alpha_1$. 
\vskip0,1truecm
Now, we have 
\begin{equation*}
\left\Vert \Phi _{n}\right\Vert _{\mathcal{H}}^{2} \geq \rho_2\left\Vert {\tilde \psi}_{n}\right\Vert _{L^{2}\left( 0,1\right) }^{2} =\rho_2\left(\alpha_2 \lambda_n\right)^{2} \int_{0}^{1} \cos^2\,\left(Nx\right)\,dx=\frac{\rho_2}{2}\left(\alpha_2 \lambda_n\right)^{2} \int_{0}^{1} \left[1+\cos\,\left(2Nx\right)\right]\,dx=\frac{\rho_2}{2}\left(\alpha_2 \lambda_n\right)^{2},
\end{equation*}
consequently, $\left( \ref{fn69}\right)$ holds.
\vskip0,1truecm
Finally, there exist sequences $\left(F_n\right)_n \subset \mathcal{H}$, $\left(\Phi _{n} \right)_n\subset D(\mathcal{A})$ and $\left(\lambda_n \right)_n\subset \mathbb{R}$ satisfying $\left( \ref{Fn}\right)$ and $\left( \ref{eq_4}\right)$. Hence, Theorem \ref{Theorem 2.1} implies that system $\left(\ref{syst_2}\right)$ is not exponentially stable.
\end{proof}

\section{Polynomial stability of $(\ref{syst1})-(\ref{cdt_1})$}

In this section, we prove that the semigroup associated to system $\left( \ref{syst_2}\right)$ is polynomially stable. Our second main result is stated as follow:
\begin{theorem}\label{Theorem 4.2} Assume that $l$ satisfies $(\ref{pi})$ and
\begin{equation}
l^2 \ne\dfrac{k_0 \rho_2 -b\rho_1}{k_0\rho_2} \left(m\pi\right)^2 - \dfrac{k \rho_1}{\rho_2 \left(k+k_0\right)} ,\quad\forall m\in  \mathbb{Z} . \label{eq_3_144}
\end{equation}
Then, for any $m\in\mathbb{N}^*$, there exists a constant $c_{m}>0$ such that
\begin{equation}
\forall \Phi _{0}\in D\left( \mathcal{A}^m\right) ,\,\,\forall t\geq 2,\,\,\left\Vert
e^{t\mathcal{A}}\Phi _{0}\right\Vert _{\mathcal{H}}\leq c_{m}\left\Vert \Phi _{0}\right\Vert
_{D\left( \mathcal{A}^m\right) }\left( \frac{\ln t}{t}\right) ^{\dfrac{m}{3}}\ln t.  \label{eq_4_1}
\end{equation}
\end{theorem}
\begin{proof} Using Theorem \ref{Theorem 2.2}, it is sufficient to show that
\begin{equation}
i\,\mathrm{I\hskip-2ptR}\subset \rho \left( \mathcal{A}\right) \label{eq_4_2}
\end{equation}
and 
\begin{equation}
\sup_{\left\vert \lambda \right\vert \,\geq \,1}\frac{1}{\lambda ^{3}}\left\Vert \left(
i\lambda I-\mathcal{A}\right) ^{-1}\right\Vert_{\mathcal{H}} <\infty .\label{eq_4_3}
\end{equation}
We start by proving $\left( \ref{eq_4_2}\right)$. Notice that, according to the fact that $0\in $ $\rho \left( \mathcal{A}\right)$ (section 2), $\mathcal{A}^{-1}$ is bounded and it is a bijection between $\mathcal{H}$ and $D(\mathcal{A})$. Since $D(\mathcal{A})$ has a compact embedding into $\mathcal{H}$, so, it follows that $\mathcal{A}^{-1}$ is a compact operator, which implies that the spectrum of $\mathcal{A}$ is discrete.
\vskip0,1truecm
Let $\lambda\in\mathbb{R}^*$ and 
\begin{equation*}
\Phi = \left(\varphi ,{\tilde{\varphi}} ,\psi ,{\tilde{\psi}} ,w ,{\tilde{w}} \right)^T\in D(\mathcal{A}) .
\end{equation*}
We prove that $i\lambda$ is not an eigenvalue of $\mathcal{A}$ by proving that the equation
\begin{equation}
\mathcal{A}\,\Phi =i\,\lambda\,\Phi  \label{eq_3_4}
\end{equation}
has only $\Phi =0$ as a solution. Assume that $\left( \ref{eq_3_4}\right)$ is true, then we have
\begin{equation}
\left\{
\begin{array}{l}
{\tilde{\varphi}}=i\lambda\varphi,\quad {\tilde{\psi}}=i\lambda\psi, \quad {\tilde{w}}=i\lambda w, \vspace{0.2cm}\\
\dfrac{k}{\rho _{1}}\left( \varphi _{x}+\psi +l\,w\right) _{x}+\dfrac{lk_{0}}{\rho _{1}}\left( w_{x}-l\varphi \right) -\dfrac{\delta}{\rho _{1}}{\tilde \varphi} =i\lambda {\tilde{\varphi}} ,\vspace{0.2cm}\\
\dfrac{b}{\rho _{2}}\psi _{xx}-\dfrac{k}{\rho _{2}}\left( \varphi _{x}+\psi
+l\,w\right) =i\lambda {\tilde{\psi}}, \vspace{0.2cm}\\
\dfrac{k_{0}}{\rho _{1}}\left( w_{x}-l\varphi \right) _{x}-\dfrac{lk}{\rho
_{1}}\left( \varphi _{x}+\psi +l\,w\right) =i\lambda {\tilde{w}} .
\end{array}
\right. \label{eq_3_5}
\end{equation}
Using $\left( \ref{dissp}\right)$, we get
\begin{equation*}
0=Re\,i\lambda \left\Vert \Phi\right\Vert_{\mathcal{H}}^2 =Re\,\left\langle i\lambda \Phi,\Phi \right\rangle _{\mathcal{H}} =Re\,\left\langle \mathcal{A}\Phi ,\Phi \right\rangle _{\mathcal{H}}= -\delta
\left\Vert \tilde{\varphi}\right\Vert _{L^{2}\left( 0,1\right) }^{2} .
\end{equation*}
So, by the first equation in $\left( \ref{eq_3_5}\right)$, we find
\begin{equation}
\varphi = {\tilde{\varphi}}=0.  \label{eq_3_10}
\end{equation}
Using $\left( \ref{eq_3_10}\right)$, we see that $\left( \ref{eq_3_5}\right)$ leads to
\begin{equation}
\left\{
\begin{array}{l}
{\tilde{\psi}}=i\lambda\psi, \quad {\tilde{w}}=i\lambda w, \vspace{0.2cm}\\
k\psi _{x}+l\left( k+k_{0}\right) w_{x}=0, \vspace{0.2cm}\\
b\psi _{xx}-k\left( \psi +l\,w\right) =-\rho _{2}\lambda^{2}\psi ,\vspace{0.2cm}\\
k_{0}w_{xx}-lk\left( \psi +l\,w\right) =-\rho _{1}\lambda^{2}w.
\end{array}
\right.  \label{eq_3_11}
\end{equation}
The third equation in $\left( \ref{eq_3_11}\right)$ implies that $k\psi +l\left( k+k_{0}\right) w$ is a constant, then, thanks to the definition of $L_*^{2}\left( 0,1\right)$, we get 
\begin{equation}
\psi =-l\left( 1+\dfrac{k_{0}}{k}\right) w. \label{eq_3_12}
\end{equation}
Using the last two equations in $\left( \ref{eq_3_11}\right)$, we obtain
\begin{equation}
lb\psi _{xx}-k_{0}w_{xx}=-\rho _{2} l\lambda^{2} \psi +\rho _{1}\lambda^{2}w.  \label{eq_3_13}
\end{equation}
Then, combining with $\left( \ref{eq_3_12}\right)$, we find
\begin{equation*}
w_{xx}+\alpha^2\lambda^{2}w=0,
\end{equation*}
where
\begin{equation}
\alpha=\sqrt{\dfrac{\rho _{2} l^2\left( k+k_{0}\right) +k\rho _{1} } { bl^2\left( k+k_{0}\right) +kk_{0} } }.\label{eq_3_130}
\end{equation}
This implies that, for $c_1 ,\,c_2\in \mathbb{C}$,
\begin{equation*}
w(x)=c_1\cos\,\left(\alpha\lambda x\right) + c_2\sin\,\left(\alpha\lambda x\right).
\end{equation*}
The boundary condition $w_{x}\left( 0\right) =0$ leads to $c_2 =0$, and then, using $\left( \ref{eq_3_12}\right)$,
\begin{equation}
\psi (x)=-l\left( 1+\dfrac{k_{0}}{k}\right)c_1\cos\,\left(\alpha\lambda x\right)\quad\hbox{and}\quad w(x)=c_1\cos\,\left(\alpha\lambda x\right).  \label{allam+}
\end{equation}
Because $\psi_x (1)=w_x (1)=0$, we have
\begin{equation*}
c_1 =0\quad\hbox{or}\quad\exists m\in \mathbb{Z} :\,\alpha\lambda =m\pi.
\end{equation*}
Assume by contradiction that 
\begin{equation}
\exists m\in \mathbb{Z} :\,\alpha\lambda =m\pi. \label{eq_3_141}
\end{equation}
Therefore, using $\left( \ref{eq_3_130}\right)$ and $\left( \ref{allam+}\right)$, we get that the last two equations in $\left( \ref{eq_3_11}\right)$ are equivalent to
\begin{equation}
\left(k_0 \rho_2 -b\rho_1\right)\lambda^2 =\dfrac{k_0}{k+k_0}\left[bl^2 \left(k+k_0\right)+kk_0\right]. \label{eq_3_142}
\end{equation}
So, combining $\left( \ref{eq_3_130}\right)$, $\left( \ref{eq_3_141}\right)$ and $\left( \ref{eq_3_142}\right)$, we get
\begin{equation*}
\exists m\in \mathbb{Z}:\,\,l^2 =\dfrac{k_0 \rho_2 -b\rho_1}{k_0\rho_2} \left(m\pi\right)^2 - \dfrac{k \rho_1}{\rho_2 \left(k+k_0\right)} , 
\end{equation*} 
which is a contraduction with $\left( \ref{eq_3_144}\right)$. Consequentely, $c_1 =0$ and hence
\begin{equation}
\psi =w=0. \label{eq_3_14}
\end{equation}
Using $\left( \ref{eq_3_14}\right)$ and the first two equations in $\left( \ref{eq_3_11}\right)$, we obtain
\begin{equation*}
{\tilde{\psi}}={\tilde{w}}=0.
\end{equation*}
Finally, $\Phi =0$ and thus
\begin{equation}
i\lambda\in \rho \left( \mathcal{A}\right). \label{eq_3_140}
\end{equation}
This ends the proof of $\left( \ref{eq_4_2}\right)$. 
\vskip0,1truecm
Now, we establish $\left( \ref{eq_4_3}\right)$ by contradiction. Assume that $\left( \ref{eq_4_3}\right)$ is false, then there exist sequences $\left(\Phi _{n}\right)_n\subset \,D\left( \mathcal{A}\right)$ and
$\left(\lambda _{n}\right)_n\subset \,\mathbb{R}$ satisfying
\begin{equation}
\left\Vert \,\Phi _{n}\right\Vert _{\mathcal{H}}\,=\,1,\quad\forall \,n\in \mathbb{N},  \label{eq_4_4}
\end{equation}%
\begin{equation}
\lim_{n\rightarrow \infty }\left\vert \lambda _{n}\right\vert =\infty \label{eq_4_5}
\end{equation}
and
\begin{equation}
\lim_{n\rightarrow \infty }\lambda _{n}^{3}\left\Vert \left( i\lambda_{n}\,I\,-\,\mathcal{A}\right) \,\Phi _{n}\right\Vert _{\mathcal{H}}\,=\,0.\label{eq_4_6}
\end{equation}
Let $\Phi _{n} =\left(\varphi_{n} ,\overset{\sim }{\varphi}_n ,\psi_n ,\overset{\sim }{\psi}_{n} ,w_n, \overset{\sim }{w}_{n} \right)^T$. The limit $\left( \ref{eq_4_6}\right)$ implies that
\begin{equation}
\left\{
\begin{array}{ll}
\lambda _{n}^{3}\left[ i\lambda _{n}\varphi _{n}-\overset{\sim }{\varphi }%
_{n}\right] \rightarrow 0\,\,&\text{in}\,\,H_{0}^{1}\left( 0,1\right),\vspace{0.2cm}\\
\lambda _{n}^{3}\left[ i\lambda _{n}\rho _{1}\overset{\sim }{\varphi }
_{n}-k\left( \varphi _{nx}+\psi _{n}+lw_{n}\right) _{x}-lk_{0}\left(
w_{nx}-l\varphi _{n}\right) +\delta {\tilde\varphi}_n \right] \rightarrow 0\,\,&\text{in}\,\,L^{2}\left( 0,1\right),\vspace{0.2cm}\\
\lambda _{n}^{3}\left[ i\lambda _{n}\psi _{n}-\overset{\sim }{\psi }_{n}
\right] \rightarrow 0\,\,&\text{in}\,\,{H_{\ast }^{1}}\left(
0,1\right),\vspace{0.2cm}\\
\lambda _{n}^{3}\left[ i\lambda _{n}\rho _{2}\overset{\sim }{\psi }
_{n}-b\psi _{nxx}+k\left( \varphi _{nx}+\psi _{n}+lw_{n}\right) \right]
\rightarrow 0\,\,&\text{in}\,\,L_*^{2}\left( 0,1\right),\vspace{0.2cm}\\
\lambda _{n}^{3}\left[ i\lambda _{n}w_{n}-\overset{\sim }{w}_{n}\right]
\rightarrow 0\,\,&\text{in}\,\,{H_{\ast }^{1}}\left( 0,1\right),\vspace{0.2cm}\\
\lambda _{n}^{3}\left[ i\lambda _{n}\rho _{1}\overset{\sim }{w}
_{n}-k_{0}\left( w_{nx}-l\varphi _{n}\right) _{x}+lk\left( \varphi
_{nx}+\psi _{n}+lw_{n}\right) \right] \rightarrow 0\,\,&\text{in}
\,\,L_*^{2}\left( 0,1\right).
\end{array}
\right.  \label{eq_4_7}
\end{equation}
We will prove that $\left\Vert \,\Phi _{n}\right\Vert_{\mathcal{H}}\,\rightarrow 0$, which gives a contradiction with 
$\left( \ref{eq_4_4}\right)$. We use several multipliers, some of them were used in \cite{afas}.
\vskip0,2truecm
\textbf{Step 1.} Using $\left( \ref{dissp}\right)$, we get
\begin{equation*}
\begin{array}{lll}
Re\left\langle \lambda _{n}^{2}\left( i\,\lambda _{n}\,I\,-\,\mathcal{A}\right) \,\,\Phi _{n},\Phi _{n}\right\rangle_{\mathcal{H}} &=& Re\left( i\lambda _{n}^{3}\left\Vert \Phi _{n}\right\Vert_{\mathcal{H}}^{2} -\lambda _{n}^{2}\left\langle \mathcal{A} \,\Phi _{n},\Phi _{n}\right\rangle_{\mathcal{H}}\right) \vspace{0.2cm}\\
&=& \delta\lambda _{n}^{2}\Vert {\tilde\varphi}_n\Vert_{L^{2}\left( 0,1\right) }^{2} .
\end{array}
\end{equation*}
So, $\left( \ref{eq_4_4}\right)$, $\left( \ref{eq_4_5}\right)$ and $\left( \ref{eq_4_6}\right)$ imply that
\begin{equation}
\lambda _{n} {\tilde\varphi}_n\longrightarrow 0\,\,\text{ in}\,\,L^{2}\left( 0,1\right). \label{eq_4_8}
\end{equation}
\vskip0,1truecm
\textbf{Step 2.} Multiplying $\left( \ref{eq_4_7}\right)_1$ by $\dfrac{1}{\lambda _{n}^2}$, and using $\left( \ref{eq_4_5}\right)$ and $\left( \ref{eq_4_8}\right)$, we obtain
\begin{equation}
\lambda _{n}^2 \varphi_{n}\longrightarrow 0\text{ in}\,\,L^{2}\left( 0,1\right). \label{eq_4_11}
\end{equation}
\vskip0,1truecm
\textbf{Step 3.} Using an integration by parts, $\left( \ref{eq_4_5}\right)$ and $\left( \ref{eq_4_7}\right)_2$, we see that
\begin{equation*}
\left\langle \left[ i\lambda _{n}\rho _{1}\overset{\sim }{\varphi }
_{n}-k\psi_{nx} -l(k+k_{0} )w_{nx} +l^2 k_0 \varphi_n +\delta {\tilde\varphi}_n\right],\lambda _{n}^2\varphi_{n}\right\rangle _{L^{2}\left( 0,1\right)}+k\lambda _{n}^2\Vert\varphi_{nx}\Vert_{L^{2}\left( 0,1\right)}^2 \longrightarrow 0,
\end{equation*}
so, using $\left( \ref{eq_4_4}\right)$, $\left( \ref{eq_4_8}\right)$ and $\left( \ref{eq_4_11}\right)$, 
\begin{equation}
\lambda _{n}\varphi_{nx}\longrightarrow 0\,\,\text{in}\,\,L^{2}\left(0,1\right). \label{eq_4_22}
\end{equation}
Moreover, because $\varphi_{n}\in H_0^1 \left(0,1\right)$, then
\begin{equation}
\lambda _{n}\varphi _{n}\longrightarrow 0\;\text{in}\,\;L^{2}\left(0,1\right),  \label{eq_4_16}
\end{equation}
and by $\left( \ref{eq_4_7}\right)_1$ and $\left( \ref{eq_4_22}\right)$, we find 
\begin{equation}
\overset{\sim }{\varphi }_{nx}\rightarrow 0\text{ in}\,\;L^{2}\left(0,1\right).  \label{eq_4_21}
\end{equation} 
\vskip0,1truecm
\textbf{Step 4.} Multiplying $\left( \ref{eq_4_7}\right)_3$ and $\left( \ref{eq_4_7}\right)_5$ by $\dfrac{1}{\lambda _{n}^{4}}$, and using $\left( \ref{eq_4_4}\right)$ and $\left( \ref{eq_4_5}\right)$, we obtain
\begin{equation}
\psi _{n}\longrightarrow 0\text{ in}\,\,L^{2}\left( 0,1\right)\quad\hbox{and}\quad w_{n}\longrightarrow 0\text{ in}\,\,L^{2}\left( 0,1\right).\label{eq_4_11*}
\end{equation}
\vskip0,1truecm
\textbf{Step 5.} Taking the inner product of $\left( \ref{eq_4_7}\right)_2$
with $\dfrac{1}{\lambda _{n}^{3}}\left[ k\psi _{nx}+l\left( k+k_{0}\right) w_{nx}\right]$ in $L^{2}\left( 0,1\right)$ and using $\left( \ref{eq_4_5}\right)$, we get
\begin{equation}
\rho _{1}\left\langle (i\rho _{1}\lambda _{n} +\delta)\overset{\sim }{\varphi }_{n},\left[
k\psi _{nx}+l\left( k+k_{0}\right) w_{nx}\right] \right\rangle_{L^{2}\left(0,1\right)}
-k\left\langle \varphi _{nxx},\left[ k\psi _{nx}+l\left( k+k_{0}\right)w_{nx}\right] \right\rangle_{L^{2}\left(0,1\right)}
\label{eq_4_28}
\end{equation}
\begin{equation*}
-\left\Vert k\psi _{nx}+l\left( k+k_{0}\right) w_{nx}\right\Vert_{L^{2}\left(0,1\right)}^{2}
+l^{2}k_{0}\left\langle \varphi _{n},\left[ k\psi _{nx}+l\left(k+k_{0}\right) w_{nx}\right] \right\rangle_{L^{2}\left(0,1\right)}
\rightarrow 0. 
\end{equation*}
Integrating by parts and using the boundary conditions, we get
\begin{equation}
\left\langle \varphi _{nxx},\left[ k\psi _{nx}+l\left( k+k_{0}\right)w_{nx}\right] \right\rangle_{L^{2}\left(0,1\right)} =-\left\langle \lambda _{n}\varphi _{nx},\left[ k\frac{\psi _{nxx}}{\lambda _{n}}+l\left( k+k_{0}\right) \frac{w_{nxx}}{\lambda _{n}}\right]\right\rangle_{L^{2}\left(0,1\right)} .   \label{eq_4_29}
\end{equation}
On the other hand, multiplying $\left( \ref{eq_4_7}\right)_4$ and $\left( \ref{eq_4_7}\right)_6$ by $\dfrac{1}{\lambda_{n}^{4}}$ and using $\left(\ref{eq_4_5}\right)$, we obtain
\begin{equation*}
\left\{
\begin{array}{ll}
i\rho _{2}\overset{\sim }{\psi }_{n}-b\dfrac{\psi _{nxx}}{\lambda_{n}}+\dfrac{k}{\lambda _{n}} \left( \varphi _{nx}+\psi _{n}+lw_{n}\right)\rightarrow 0\,&\text{\ in}\,\, L^{2}\left( 0,1\right), \vspace{0.2cm}\\
i\rho _{1}\overset{\sim }{w}_{n}-k_{0}\dfrac{w_{nxx}}{\lambda _{n}}+lk_{0}\dfrac{\varphi _{nx}}{\lambda _{n}}+\dfrac{lk}{\lambda _{n}}\left(\varphi _{nx}+\psi _{n}+lw_{n}\right) \rightarrow 0\,&\text{ in}\,\,L^{2}\left( 0,1\right).
\end{array}
\right.
\end{equation*}
Exploiting $\left( \ref{eq_4_4}\right)$, we get
\begin{equation}
\left( \dfrac{1}{\lambda _{n}}\psi _{nxx}\right)_{n}\,\,\hbox{and} \,\, \left(\dfrac{1}{\lambda _{n}} w_{nxx}\right)_{n}\,\,\text{are bounded in}\,\,L^{2}\left( 0,1\right),\label{eq_4_30}
\end{equation}
then, using  $\left( \ref{eq_4_8}\right)$, $\left( \ref{eq_4_22}\right)$, $\left( \ref{eq_4_29}\right)$ and $\left( \ref{eq_4_30}\right)$, we deduce that
\begin{equation}
\left\langle \varphi _{nxx},\left[ k\psi _{nx}+l\left( k+k_{0}\right)
w_{nx}\right] \right\rangle_{L^{2}\left(0,1\right)} \rightarrow 0,  \label{eq_4_31}
\end{equation}
so, exploiting $\left( \ref{eq_4_4}\right)$, 
$\left(\ref{eq_4_16}\right)$, $\left( \ref{eq_4_8}\right)$ and $\left( \ref{eq_4_28}\right)$, we have 
\begin{equation}
k\psi _{nx}+l\left( k+k_{0}\right) w_{nx}\rightarrow 0\text{ in}\,\, L^{2}\left( 0,1\right). \label{eq_4_32}
\end{equation}
\vskip0,1truecm
\textbf{Step 6.} Taking the inner product of $\left( \ref{eq_4_7}\right)_4$ with $\dfrac{\psi _{n}}{\lambda _{n}^{3}}$ in $L^{2}\left( 0,1\right)$, using $\left( \ref{eq_4_5}\right)$, integrating by parts and using the boundary conditions, we obtain
\begin{equation*}
-\rho _{2}\left\langle \overset{\sim }{\psi }_{n},\left( i\lambda _{n}\psi
_{n}-\overset{\sim }{\psi }_{n}\right) \right\rangle_{L^{2}\left(0,1\right)} -\rho _{2}\left\Vert
\overset{\sim }{\psi }_{n}\right\Vert_{L^{2}\left(0,1\right)} ^{2}
+b\left\Vert \psi_{nx}\right\Vert_{L^{2}\left(0,1\right)} ^{2}+k\left\langle \left( \varphi _{nx}+\psi_{n}+lw_{n}\right) ,\psi _{n}\right\rangle_{L^{2}\left(0,1\right)} \rightarrow 0,
\end{equation*}
then, using $\left( \ref{eq_4_4}\right)$, $\left( \ref{eq_4_5}\right)$, $\left( \ref{eq_4_7}\right)_3$ and 
$\left( \ref{eq_4_11*}\right)$, we find
\begin{equation}
b\left\Vert \psi _{nx}\right\Vert_{L^{2}\left(0,1\right)}^{2}-\rho _{2}\left\Vert \overset{\sim }{%
\psi }_{n}\right\Vert_{L^{2}\left(0,1\right)}^{2}\rightarrow 0 .  \label{eq_4_33}
\end{equation}
On the other hand, taking the inner product of $\left( \ref{eq_4_7}\right)_6$ with $\dfrac{w_{n}}{\lambda _{n}^{3}}$ in 
$L^{2}\left( 0,1\right)$, using $\left( \ref{eq_4_5}\right)$,
integrating by parts and using the boundary conditions, we observe that
\begin{equation*}
-\rho _{1}\left\langle \overset{\sim }{w}_{n},\left( i\lambda _{n}w_{n}-\overset{\sim }{w}_{n}\right)\right\rangle_{L^{2}\left(0,1\right)} -\rho _{1}\left\Vert \overset{
\sim }{w}_{n}\right\Vert_{L^{2}\left(0,1\right)}^{2} +k_{0}\left\Vert w_{nx}\right\Vert_{L^{2}\left(0,1\right)}^{2}
\end{equation*}
\begin{equation*}
+lk_{0}\left\langle \varphi_{nx},w_{n}\right\rangle_{L^{2}\left(0,1\right)} +lk\left\langle \left( \varphi _{nx}+\psi
_{n}+lw_{n}\right) ,w_{n}\right\rangle_{L^{2}\left(0,1\right)} \rightarrow 0.
\end{equation*}
By $\left( \ref{eq_4_4}\right)$, $\left( \ref{eq_4_5}\right)$, $\left( \ref{eq_4_7}\right)_5$ and $\left( \ref{eq_4_11*}\right)$, we deduce that
\begin{equation}
k_{0}\left\Vert w_{nx}\right\Vert_{L^{2}\left(0,1\right)}^{2}-\rho _{1}\left\Vert \overset{\sim }{w%
}_{n}\right\Vert_{L^{2}\left(0,1\right)}^{2}\rightarrow 0.  \label{eq_4_34}
\end{equation}
\vskip0,1truecm
\textbf{Step 7.} Taking the inner product of $\left( \ref{eq_4_7}\right)_4$ with $\dfrac{w_{n}}{\lambda _{n}^{3}}$ and of 
$\left( \ref{eq_4_7}\right)_6$ with $\dfrac{\psi _{n}}{\lambda _{n}^{3}}$ in $L^{2}\left( 0,1\right)$, and using $\left( \ref{eq_4_5}\right)$, we get
\begin{equation*}
\left\{
\begin{array}{l}
\left\langle \left[ i\lambda _{n}\rho _{2}\overset{\sim }{\psi }_{n}-b\psi
_{nxx}+k\left( \varphi _{nx}+\psi _{n}+lw_{n}\right) \right]
,w_{n}\right\rangle_{L^{2}\left(0,1\right)} \rightarrow 0,  \vspace{0.2cm}\\
\left\langle \left[ i\lambda _{n}\rho _{1}\overset{\sim }{w}_{n}-k_{0}\left(
w_{nx}-l\varphi _{n}\right) _{x}+lk\left( \varphi _{nx}+\psi
_{n}+lw_{n}\right) \right] ,\psi _{n}\right\rangle_{L^{2}\left(0,1\right)} \rightarrow 0.
\end{array}
\right.
\end{equation*}
Integrating by parts and using the boundary conditions, we obtain
\begin{equation*}
-\rho _{2}\left\langle \overset{\sim }{\psi }_{n},\left( i\lambda _{n}w_{n}-%
\overset{\sim }{w}_{n}\right) \right\rangle_{L^{2}\left(0,1\right)} -\rho _{2}\left\langle \overset{%
\sim }{\psi }_{n},\overset{\sim }{w}_{n}\right\rangle_{L^{2}\left(0,1\right)}
\end{equation*}
\begin{equation*}
+b\left\langle \psi_{nx},w_{nx}\right\rangle_{L^{2}\left(0,1\right)} +k\left\langle \left( \varphi _{nx}+\psi
_{n}+lw_{n}\right) ,w_{n}\right\rangle_{L^{2}\left(0,1\right)} \rightarrow 0
\end{equation*}
and
\begin{equation*}
-\rho _{1}\left\langle \overset{\sim }{w}_{n},\left( i\lambda _{n}\psi _{n}-
\overset{\sim }{\psi }_{n}\right) \right\rangle_{L^{2}\left(0,1\right)} -\rho _{1}\left\langle
\overset{\sim }{w}_{n},\overset{\sim }{\psi }_{n}\right\rangle_{L^{2}\left(0,1\right)}
-lk_{0}\left\langle \varphi _{n} ,\psi_{nx}\right\rangle_{L^{2}\left(0,1\right)}
\end{equation*}
\begin{equation*}
+k_{0}\left\langle w_{nx} ,\psi_{nx}\right\rangle_{L^{2}\left(0,1\right)} +lk\left\langle \left( \varphi _{nx}+\psi
_{n}+lw_{n}\right) ,\psi _{n}\right\rangle_{L^{2}\left(0,1\right)} \rightarrow 0,
\end{equation*}
then, using $\left( \ref{eq_4_4}\right)$, $\left( \ref{eq_4_5}\right)$, $\left( \ref{eq_4_7}\right)_3$, 
$\left( \ref{eq_4_7}\right)_5$, $\left( \ref{eq_4_11}\right)$ and $\left( \ref{eq_4_11*}\right)$, we obtain
\begin{equation*}
\left\{
\begin{array}{l}
-\rho _{2}\left\langle \overset{\sim }{\psi }_{n},\overset{\sim }{w}_{n}\right\rangle_{L^{2}\left(0,1\right)} +b\left\langle \psi _{nx},w_{nx}\right\rangle_{L^{2}\left(0,1\right)}\rightarrow 0,\vspace{0.2cm}\\
-\rho _{1}\left\langle \overset{\sim }{\psi }_{n} ,\overset{\sim }{w}_{n}\right\rangle_{L^{2}\left(0,1\right)} +k_{0}\left\langle \psi _{nx} ,w_{nx}\right\rangle_{L^{2}\left(0,1\right)}\rightarrow 0,
\end{array}
\right.
\end{equation*}
which implies that
\begin{equation}
\left( \frac{\rho _{2}}{b}-\frac{\rho _{1}}{k_{0}}\right)
\left\langle \overset{\sim }{\psi }_{n},\overset{\sim }{w}_{n}\right\rangle_{L^{2}\left(0,1\right)} \rightarrow 0  \label{eq_4_35}
\end{equation}
and
\begin{equation}
\left( \frac{b}{\rho _{2}}-\frac{k_{0}}{\rho _{1}}\right)\left\langle \psi _{nx},w_{nx}\right\rangle_{L^{2}\left(0,1\right)} \rightarrow 0. \label{eq_4_36}
\end{equation}
\vskip0,1truecm
\textbf{Step 8.} We distinguish in this step two cases.
\vskip0,1truecm
\textbf{Case 1:} $\dfrac{b}{\rho _{2}} \ne \dfrac{k_{0}}{\rho _{1}}$. From $\left( \ref{eq_4_35}\right)$ and $\left( \ref{eq_4_36}\right)$, we see that
\begin{equation}
\left\langle \overset{\sim }{\psi }_{n},\overset{\sim }{w}_{n}\right\rangle_{L^{2}\left(0,1\right)} \rightarrow 0 \quad\hbox{and}\quad\left\langle \psi _{nx},w_{nx}\right\rangle_{L^{2}\left(0,1\right)} \rightarrow 0. \label{eq_4_360}
\end{equation}
Therefore, taking the inner product in $L^{2}\left( 0,1\right)$ of $\left( \ref{eq_4_32}\right)$, first, with $\psi_{nx}$, and second, with $w_{nx}$, we obtain
\begin{equation}
\psi _{nx}\rightarrow 0\quad\hbox{and}\quad w_{nx}\rightarrow 0 \text{ in}\,\, L^{2}\left( 0,1\right),\label{eq_4_37}
\end{equation}
and then, by $\left( \ref{eq_4_33}\right)$, $\left( \ref{eq_4_34}\right)$ and $\left( \ref{eq_4_37}\right)$,
\begin{equation}
\overset{\sim }{\psi }_{n}\rightarrow 0\quad\hbox{and}\quad
\overset{\sim }{w}_{n}\rightarrow 0\text{ in}\,\, L^{2}\left( 0,1\right).\label{eq_4_38}
\end{equation}
Finally, combining $\left( \ref{eq_4_8}\right)$, $\left( \ref{eq_4_22}\right)$, $\left( \ref{eq_4_16}\right)$, $\left( \ref{eq_4_11*}\right)$, $\left( \ref{eq_4_37}\right)$ and $\left( \ref{eq_4_38}\right)$, we get
\begin{equation}
\left\Vert \Phi _{n}\right\Vert _{\mathcal{H}}\longrightarrow 0, \label{eq_4_380}
\end{equation}
which is a contradiction with $\left( \ref{eq_4_4}\right)$, so $\left( \ref{eq_4_3}\right)$ holds. Consequentely,
$\left( \ref{eq_4_1}\right)$ is satisfied.
\vskip0,1truecm
\textbf{Case 2:} $\dfrac{b}{\rho _{2}}=\dfrac{k_{0}}{\rho _{1}}$. Using $\left( \ref{eq_4_7}\right)_4$ and $\left( \ref{eq_4_7}\right)_6$, we obtain 
\begin{equation*}
\left\{
\begin{array}{ll}
\lambda _{n}^{2}\left[ -\dfrac{i\rho _{2}}{b}\lambda _{n} \left(i\lambda _{n} \psi _{n} -{\tilde \psi}_n\right) -\dfrac{\rho _{2}}{b}\lambda _{n}^{2}\psi _{n}-\psi
_{nxx}+\dfrac{k}{b}\left( \varphi _{nx}+\psi _{n}+lw_{n}\right) \right]
\rightarrow 0\,\,&\text{in}\,\,L^{2}\left( 0,1\right), \vspace{0.2cm}\\
\lambda _{n}^{2}\left[-\dfrac{i\rho _{2}}{b}\lambda _{n} \left(i\lambda _{n} w_{n} -{\tilde w}_n\right) -\dfrac{\rho _{2}}{b}\lambda _{n}^{2}w_{n}-\left(
w_{nx}-l\varphi _{n}\right) _{x}+\dfrac{lk}{k_{0}}\left( \varphi
_{nx}+\psi _{n}+lw_{n}\right) \right] \rightarrow 0\,\,&\text{ in}\,\,L^{2}\left( 0,1\right),
\end{array}
\right.  
\end{equation*} 
so, using $\left( \ref{eq_4_7}\right)_3$ and $\left( \ref{eq_4_7}\right)_5$, we find
\begin{equation}
\left\{
\begin{array}{ll}
\lambda _{n}^{2}\left[ -\dfrac{\rho _{2}}{b}\lambda _{n}^{2}\psi _{n}-\psi
_{nxx}+\dfrac{k}{b}\left( \varphi _{nx}+\psi _{n}+lw_{n}\right) \right]
\rightarrow 0\,\,&\text{in}\,\,L^{2}\left( 0,1\right), \vspace{0.2cm}\\
\lambda _{n}^{2}\left[ -\dfrac{\rho _{2}}{b}\lambda _{n}^{2}w_{n}-\left(
w_{nx}-l\varphi _{n}\right) _{x}+\dfrac{lk}{k_{0}}\left( \varphi
_{nx}+\psi _{n}+lw_{n}\right) \right] \rightarrow 0\,\,&\text{ in}\,\,L^{2}\left( 0,1\right).
\end{array}
\right.  \label{eq_4_39}
\end{equation}
Then, using $\left( \ref{eq_4_5}\right)$, $\left( \ref{eq_4_22}\right)$ and $\left( \ref{eq_4_11*}\right)$, we get
\begin{equation}
\dfrac{\rho _{2}}{b}\lambda _{n}^{2}\psi _{n}+\psi _{nxx}\rightarrow 0\,\,\text{in}\,\,L^{2}\left( 0,1\right)\quad\hbox{and}\quad 
\dfrac{\rho _{2}}{b}\lambda _{n}^{2}w_{n}+w_{nxx}\rightarrow 0\,\,\text{in}\,\,L^{2}\left( 0,1\right).\label{eq_4_40}
\end{equation}
Multiplying $\left( \ref{eq_4_40}\right)_1$ by $k$ and $\left(\ref{eq_4_40}\right)_2$ by $l(k+k_{0})$ and adding the obtained limits, and multiplying $\left( \ref{eq_4_40}\right)_1$ by $k$ and $\left( \ref{eq_4_40}\right)_2$ by $-l(k+k_{0})$ and adding the limits, we obtain
\begin{equation}
\left\{
\begin{array}{l}
\dfrac{\rho _{2}}{b}\lambda _{n}^{2}\left[ k\psi _{n}+l(k+k_{0})w_{n}\right]
+\left[ k\psi _{nxx}+l(k+k_{0})w_{nxx}\right] \rightarrow 0\,\,\,\text{in}\,\,L^{2}\left( 0,1\right), \vspace{0.2cm}\\
\dfrac{\rho _{2}}{b}\lambda _{n}^{2}\left[ k\psi _{n}-l(k+k_{0})w_{n}\right]
+\left[ k\psi _{nxx}-l(k+k_{0})w_{nxx}\right] \rightarrow 0\,\,\text{ in}\,\,L^{2}\left( 0,1\right).
\end{array}
\right. \label{eq_4_41}
\end{equation}
Taking the inner product in $L^{2}\left( 0,1\right)$ of $\left( \ref{eq_4_41}\right)_1$ and $\left( \ref{eq_4_41}\right)_2$
with $\left[ k\psi _{n}+l(k+k_{0})w_{n}\right]$, integrating by parts and using the boundary conditions, we get
\begin{equation*}
\dfrac{\rho _{2}}{b}\left\Vert k\lambda _{n}\psi _{n}+l(k+k_{0})\lambda
_{n}w_{n}\right\Vert_{L^{2}\left(0,1\right)}^{2} -\Vert k\psi _{nx}+l(k+k_{0})w_{nx}\Vert_{L^{2}\left(0,1\right)}^2
\rightarrow 0
\end{equation*}
and
\begin{equation*}
\dfrac{\rho_{2}}{b}\left\langle \lambda_{n}^{2}\left[ k\psi
_{n}-l(k+k_{0})w_{n}\right] ,\left[ k\psi _{n}+l(k+k_{0})w_{n}\right]
\right\rangle_{L^{2}\left(0,1\right)}
\end{equation*}
\begin{equation*}
-\left\langle \left[ k\psi _{nx}-l(k+k_{0})w_{nx}\right] ,
\left[ k\psi _{nx}+l(k+k_{0})w_{nx}\right] \right\rangle_{L^{2}\left(0,1\right)} \rightarrow 0,
\end{equation*}
then, using $\left( \ref{eq_4_4}\right)$ and $\left( \ref{eq_4_32}\right)$, we obtain
\begin{equation}
\left\{
\begin{array}{l}
k\lambda _{n}\psi _{n}+l(k+k_{0})\lambda _{n}w_{n}\rightarrow 0\,\,\text{in}%
\,\,L^{2}\left( 0,1\right),\vspace{0.2cm}\\
k^{2}\left\Vert \lambda _{n}\psi _{n}\right\Vert_{L^{2}\left(0,1\right)}^{2}-l^{2}(k+k_{0})^{2}\left\Vert \lambda_{n}w_{n}
\right\Vert_{L^{2}\left(0,1\right)}^{2}\rightarrow 0.
\end{array}
\right.  \label{eq_4_42}
\end{equation}
Taking the inner product in $L^{2}\left( 0,1\right)$ of $\left( \ref{eq_4_39}\right)_1$ with $w_{n}$, and $\left( \ref{eq_4_39}\right)_2$ with $\psi _{n}$, integrating by parts and using the boundary conditions, we get
\begin{equation}
-\dfrac{\rho _{2}}{b}\lambda _{n}^{4}\left\langle \psi
_{n},w_{n}\right\rangle_{L^{2}\left(0,1\right)} +\lambda _{n}^{2} \left\langle \psi
_{nx},w_{nx}\right\rangle_{L^{2}\left(0,1\right)} -\dfrac{k}{b}\left\langle \lambda_{n}^2\varphi_{n},w_{nx}\right\rangle_{L^{2}\left(0,1\right)} \label{eq_4_43}
\end{equation}
\begin{equation*}
+\dfrac{k}{b}\left\langle \lambda _{n}\psi _{n},\lambda _{n}w_{n}\right\rangle_{L^{2}\left(0,1\right)} +\dfrac{lk%
}{b}\left\Vert \lambda _{n}w_{n}\right\Vert_{L^{2}\left(0,1\right)}^{2}\rightarrow 0
\end{equation*}
and
\begin{equation}
-\dfrac{\rho _{2}}{b}\lambda _{n}^{4}\left\langle \psi_{n} ,w_{n}\right\rangle_{L^{2}\left(0,1\right)} +\lambda _{n}^{2}
\left\langle\psi_{nx} ,w_{nx}\right\rangle_{L^{2}\left(0,1\right)} -l\left( 1+\dfrac{k}{k_{0}}\right)\left\langle
\psi _{nx} ,\lambda_{n}^2\varphi _{n}\right\rangle_{L^{2}\left(0,1\right)} \label{eq_4_430}
\end{equation}
\begin{equation*}
+\dfrac{lk}{k_{0}}\left\Vert \lambda _{n}\psi _{n}\right\Vert_{L^{2}\left(0,1\right)}^{2}+\dfrac{l^{2}k}{k_{0}}%
\left\langle\lambda _{n}\psi _{n}, \lambda _{n}w_{n}\right\rangle_{L^{2}\left(0,1\right)} \rightarrow 0,
\end{equation*}
then, multiplying $\left( \ref{eq_4_43}\right)$ by $\dfrac{bk_0}{k}$, and $\left( \ref{eq_4_430}\right)$ by $\dfrac{-bk_0}{k}$, adding the obtained limits and using $\left( \ref{eq_4_4}\right)$ and $\left( \ref{eq_4_11}\right)$, we find
\begin{equation}
lk_{0}\left\Vert \lambda _{n}w_{n}\right\Vert_{L^{2}\left(0,1\right)}^{2}-lb\left\Vert \lambda
_{n}\psi _{n}\right\Vert_{L^{2}\left(0,1\right)}^{2}+\left( k_{0}-l^{2}b\right)
\left\langle \lambda _{n}\psi_{n},\lambda _{n} w_{n}\right\rangle_{L^{2}\left(0,1\right)}
\rightarrow 0.  \label{eq_4_44}
\end{equation}
Multiplying $\left( \ref{eq_4_7}\right)_3$ and $\left( \ref{eq_4_7}\right)_5$ by $\dfrac{1}{\lambda _{n}^{2}}$, and using
$\left( \ref{eq_4_4}\right)$ and $\left( \ref{eq_4_5}\right)$, we have
\begin{equation}
\left( \lambda _{n}\psi _{n}\right)_{n} \,\,\hbox{and}\,\,
\left( \lambda _{n}w_{n}\right) _{n}\,\,\text{are bounded in}\,\,L^{2}\left( 0,1\right).\label{eq_4_25}
\end{equation}
So, by taking the inner product in $L^{2}\left( 0,1\right)$ of $\left( \ref{eq_4_42}\right)_1$ with $\lambda_n \psi_n$, and using $\left( \ref{eq_4_25}\right)$, we have  
\begin{equation}
k\left\Vert \lambda _{n}\psi _{n}\right\Vert_{L^{2}\left(0,1\right)}^{2}+l(k+k_{0})
\left\langle \lambda _{n}w_{n},\lambda _{n}\psi _{n}\right\rangle_{L^{2}\left(0,1\right)}
\rightarrow 0. \label{eq_4_45-}
\end{equation}
Combining $\left( \ref{eq_4_42}\right)_2$ and $\left( \ref{eq_4_44}\right)$, we get
\begin{equation}
\dfrac{1}{l(k+k_{0})^{2}}\left[ k_{0}k^{2}-bl^{2}(k+k_{0})^{2}\right]
\left\Vert \lambda _{n}\psi _{n}\right\Vert_{L^{2}\left(0,1\right)}^{2}+\left( k_{0}-l^{2}b\right)
\left\langle \lambda _{n}w_{n},\lambda _{n}\psi _{n}\right\rangle_{L^{2}\left(0,1\right)}
\rightarrow 0,\label{eq_4_45}
\end{equation}
so, multiplying $\left( \ref{eq_4_45-}\right)$ by $\dfrac{\left(k+k_0\right)\left(k_0 -l^2 b\right)}{k_0}$, and 
$\left( \ref{eq_4_45}\right)$ by $\dfrac{-l\left(k+k_0 \right)^2}{k_0}$, adding the obtained limits and noting that $\dfrac{b}{\rho _{2}}=\dfrac{k_{0}}{\rho _{1}}$, we obtain
\begin{equation*}
\left[ kk_{0}+bl^{2}(k+k_{0})\right] \;\left\Vert \lambda _{n}\psi
_{n}\right\Vert_{L^{2}\left(0,1\right)}^{2}\rightarrow 0.
\end{equation*}
Then
\begin{equation}
\lambda _{n}\psi_{n}\rightarrow 0\,\,\text{in}\,\, L^{2}\left( 0,1\right)\label{eq_4_46}
\end{equation}
and, using $\left( \ref{eq_4_42}\right)_1$,
\begin{equation}
\lambda _{n}w_{n}\rightarrow 0\,\,\text{in}\,\, L^{2}\left( 0,1\right). \label{eq_4_47}
\end{equation}
Using $\left( \ref{eq_4_5}\right)$, $\left( \ref{eq_4_7}\right)_3$, $\left( \ref{eq_4_7}\right)_5$,
$\left( \ref{eq_4_46}\right)$ and $\left( \ref{eq_4_47}\right)$, we deduce that
\begin{equation}
\overset{\sim }{\psi }_{n}\rightarrow 0\,\,\text{in}\,\, L^{2}\left( 0,1\right)\quad\hbox{and}\quad \overset{\sim }{w}_{n}\rightarrow 0\,\,\text{in}\,\,L^{2}\left( 0,1\right).
\label{eq_4_48}
\end{equation}
Taking the inner product in $L^{2}\left( 0,1\right)$ of $\left( \ref{eq_4_40}\right)_1$ with $\psi_{n}$, and $\left( \ref{eq_4_40}\right)_2$ with $w_{n}$, integrating by parts and using the boundary conditions, we get
\begin{equation*}
\dfrac{\rho _{2}}{b}\left\Vert \lambda _{n}\psi _{n}\right\Vert_{L^{2}\left(0,1\right)}^{2}-\left\Vert \psi _{nx}
\right\Vert_{L^{2}\left(0,1\right)}^{2}\rightarrow 0\quad\hbox{and}\quad 
\dfrac{\rho _{2}}{b}\left\Vert \lambda _{n}w_{n}\right\Vert_{L^{2}\left(0,1\right)}^{2}-\left\Vert
w_{nx}\right\Vert_{L^{2}\left(0,1\right)}^{2}\rightarrow 0,
\end{equation*}
then, from $\left( \ref{eq_4_46}\right)$ and $\left( \ref{eq_4_47}\right)$, we conclude that
\begin{equation}
\psi _{nx}\rightarrow 0\,\,\text{in}\,\, L^{2}\left( 0,1\right)\quad\hbox{and}\quad w_{nx}\rightarrow 0\,\,\text{in}\,\, L^{2}\left( 0,1\right).
\label{eq_4_49}
\end{equation}
Finally, $\left( \ref{eq_4_8}\right)$, $\left( \ref{eq_4_22}\right)$, $\left( \ref{eq_4_16}\right)$, $\left( \ref{eq_4_11*}\right)$, $\left( \ref{eq_4_48}\right)$ and $\left( \ref{eq_4_49}\right)$ imply $\left( \ref{eq_4_380}\right)$,
which is a contradiction with $\left( \ref{eq_4_4}\right)$. Consequentely, in both cases $\dfrac{b}{\rho _{2}}\ne\dfrac{k_{0}}{\rho _{1}}$ and $\dfrac{b}{\rho _{2}}=\dfrac{k_{0}}{\rho _{1}}$, $\left( \ref{eq_4_3}\right)$ holds, and so $\left( \ref{eq_4_1}\right)$ is satisfied. Hence, the proof of Theorem \ref{Theorem 4.2} is completed.
\end{proof}

\end{document}